\documentclass{article}
\usepackage{bbm}
\usepackage{amsmath,natbib} 
\def\N{\mathbbm{N}}
\def\R{\mathbbm{R}}
\def\1{\mathbbm{1}}
\def\M{\mathcal{M}}
\def\P{\mathbbm{P}}
\def\Q{\mathbbm{Q}}
\def\KT{\mathbbm{KT}}
\def\Z{\mathbbm{Z}}
\def\pen{\mathrm{pen}}
\DeclareMathOperator*{\argmax}{arg\,max}

\DeclareMathOperator*{\kl}{\operatorname{KL}}
\newtheorem{lem}{Lemma}
\newtheorem{prop}{Proposition}
\newtheorem{theo}{Theorem}

\def\pen{\mathrm{pen}}
\renewcommand{\star}{*}
\author{
Antonio Galves, Universidade de S$\tilde{a}$o Paulo\\
Aur\'elien Garivier, CNRS UMR 5141 \& Telecom ParisTech\\
Elisabeth Gassiat, Universit\'e Paris Sud CNRS UMR 8628}
\title{Joint estimation of intersecting context tree models}
\date{May 4, 2012}

\begin{document}
\maketitle

\abstract

We study a problem of model selection for data produced by two different context tree sources. 
Motivated by linguistic questions, we consider the case where the probabilistic context trees corresponding to the two sources are finite and share many of their contexts. In order to understand the differences between the two sources, it is important to identify which contexts and which transition probabilities are specific to each source. 
We consider a class of probabilistic context tree models with three types of contexts: those which appear in one, the other, or both sources.
We use a BIC penalized maximum likelihood procedure that jointly  estimates the two sources.
We propose a new algorithm which efficiently computes the estimated context trees.
We prove that the procedure is strongly consistent.
We also present a simulation study showing the practical advantage of our procedure over a procedure that works separately on each dataset.

{\bf Key words:} BIC, context tree Models, joint estimation, penalized maximum likelihood,  variable length markov chains.

\section{Introduction}

We assign probabilistic context tree models to data produced by two different sources on the same finite alphabet $A$. 
Probabilistic context tree models were first introduced in \cite{rissanen1983} as a flexible and parsimonious model for data compression. Originally  called by Rissanen {\sl   finite memory source} or {\sl probabilistic tree},
this class of models recently became popular in the statistics literature under the name of {\sl Variable Length Markov Chains   (VLMC)} \cite{buhlmann1999}.
The idea behind the notion of variable memory models is that, given the whole past,  the
conditional  distribution of each symbol  only depends on a finite part of the past and the length of this relevant portion is a function of the past itself. Following Rissanen we 
call \emph{context} the minimal relevant part of each past. The set of all contexts satisfies the suffix property which means that no context is a proper suffix of another context. This property allows us to represent the set of all contexts as a rooted labeled tree, by reading the contexts' symbols from the nodes to the root. With this representation, the process is described by the tree of all contexts, called context tree, together with a family of probability measures on $A$ indexed by the contexts. 
In this work we shall only consider finite context trees.
The probability distribution of a context gives the transition probability to the next symbol from any past having this context as a suffix. From now on, the pair composed by the context tree and its family of probability measures will be called \emph{probabilistic context tree}.

The issue we consider here was suggested by a linguistic case study presented in~\cite{galves2009}. This paper addresses the problem of characterizing rhythmic patterns displayed by two variants of Portuguese:  Brazilian  and  European. This is done by considering two data sets consisting of encoded newspaper texts in two languages. Each data set was analysed separately using a penalized maximum likelihood procedure which selected two different probabilistic context trees corresponding to the two variants of Portuguese. A striking feature emerging from this analysis is the fact that most of the contexts and corresponding transition probabilities are common to the two dialects of Portuguese. Obviously the discriminant features characterizing the different rhythms implemented by the two dialects are expressed by the contexts which appear in one but not in the other model. 

To identify those discriminant contexts, the first idea is to estimate separately the context tree for each set of observations, using some classical context tree estimator like the algorithm Context \cite{rissanen1983} or a penalized maximum likelihood procedure as in \cite{csiszar2006} (see also \cite{garivierLeonardi11CTestimation}), and then compare the obtained trees. This is precisely what is done in \cite{galves2009}.
However, such an approach does not use the information that the two sources share some identical contexts and probability distributions. We propose in this paper a selection method using penalized maximum likelihood for the whole set of observations.

In this paper, we argue that a joint model selection more efficiently identifies the relevant features and estimates the parameters. 
The joint estimation of the two probabilistic context trees is accomplished by a penalized maximum likelihood criterium. Namely, we distinguish two types of contexts: those which appear in both sources with the same probability distribution (we call them \emph{shared contexts}), and the others. The latters appear either in only one of the two sources, or appear in both sources but with different associated probability distributions. 

At first sight the huge number of models in the class suggests that such a procedure is intractable. Actually this is not the case. We show that the Context Tree Maximizing procedure, which has been described in  \cite{willems1995}, can be adapted to recursively find the maximizer: we propose a new algorithm to efficiently compute the estimated context trees. We prove the strong consistency of the procedure. Our proof is inspired by some arguments given in~\cite{csiszar2006}, which handles the case of a single (but possibly infinite) context tree source estimation; as is \cite{garivier2006}, the size of the trees is not bounded in the maximization procedure. We also present a simulation study showing the significant advantage of our procedure, for the estimation of the shared contexts, over a procedure that works separately on each dataset.

The paper is organized as follows.
In Section~\ref{sec:notations}, we present the joint context tree estimation problem and the notation.
Section~\ref{sec:estimator} is devoted to the presentation of the penalized maximum likelihood estimator we study in this paper. For an appropriate choice of the penalty function, a strong consistency result is given.
We describe in Section~\ref{sec:algo} how to efficiently compute the joint estimator. This is a challenging task, as the number of possible models grows exponentially with the sample size. We show how to take advantage of the recursive tree structure to build an efficient greedy algorithm. % that operates in linear time.
The value of this estimator is experimentally shown in Section~\ref{sec:simu} through a simulation study.
The proof of the consistency result is given in Appendix~\ref{sec:proof}. It relies on a technical result on the Krichevsky-Trofimov distribution that is given in Appendix~\ref{sec:tech}.

\section{Notation}\label{sec:notations}
Let $A$ be a finite alphabet, and $A^{\star}=\cup_{n\in\N}A^{n}$ the set of all possible strings including the empty string $\epsilon$. 
Denote also by $A^{+}=\cup_{n\geq 1}A^{n}$ the set of non-empty strings. A string $s\in A^+$ has \emph{length} $|s|=n$ if $s\in A^{n}$, and we note $s=s_{1:|s|}$. The empty string has length $0$.
The \emph{concatenation} of strings $s$ and $s'$ is denoted by $ss'$. $s'$ is a \emph{suffix} of $s$ if there exists a string $u$ such that $s=us'$; it is a \emph{proper} suffix if $u\neq \epsilon$.

A \emph{tree} $\tau$ is a non-empty  subset of $A^{\star}$ such that no $s_{1}\in\tau$ is a suffix of any other $s_{2}\in\tau$. The \emph{depth} of a finite tree $\tau$ is defined as 
\[
%\mathrm{depth}
D(\tau) = \max \Big\{ |s| \;:\; {s\in\tau} \Big\}\;.\]
A tree is \emph{complete} if each node except the leaves has exactly $|A|$ children (here $|A|$ denotes the number of elements in $A$). Note that $\{\epsilon\}$ is a complete tree. 
%(If $|A|=2$, complete trees are irreducible trees). 

Let ${\cal P}_{A}$ be the  $(|A|-1)$-dimensional simplex, that is the subset of vectors $p=(p_{a})_{a\in A}$ in $\R^{|A|}$ such that $p_{a}\geq 0$, $a\in A$ and
$\sum_{a\in A}p_{a}=1$.
\\
To define a stationary context tree source, we need a complete tree $\tau$ and a parameter $\theta\in {\cal P}_{A}^{\tau}$,  that is $\theta=(\theta (s))_{s\in\tau}$ where, for any $s\in\tau$,  $\theta (s)\in{\cal P}_{A}$.
The $A$-valued stochastic process $Z=(Z_{n})_{n\in\Z}$ is said to be a stationary context-tree source (or variable length Markov Chain) with distribution $\P_{\tau,\theta}$ if for any semi-infinite sequence denoted by $z_{-\infty:0}$, there exists one (and only one) $s\in \tau$ such that $s$ is a suffix of $z_{-\infty:-1}$, and such that, for any $n\geq |s|$, if the event  $\{ Z_{-n:-1}=z_{-n:-1}\}$ has positive probability, the conditional distribution of $Z_{0}$ given 
$\{ Z_{-n:-1}=z_{-n:-1}\}$ is $\theta (s)$ and thus depends only on $z_{-|s|:-1}$.
%and moreover no proper suffix of $s$ has this property. 
Following Rissanen, an element of $\tau$ is called a \emph{context}.
In the case when $\tau =\{\epsilon\}$, the source is called \emph{memoryless}.

For  any $s\in \tau$, any integer $n$ and any $z_{1:n}\in A^n$, denote by $S(s;z_{1:n})$ the string with the symbols that appear after an occurrence of $s$ in the sequence $z_{1:n}$. Formally,
\[S(s; z_{1:n}) = \bigodot_{i: z_{i-|s|:i-1}=s} z_i \;,\]
where $\odot$ denotes the concatenation operator. When $z_{i-|s|:i-1}=s$, we say that $z_{1}$ is in context $s$.
Besides, denote by $I(z_{1:n};\tau)$ the set of indices $i$ of $z_{1:n}$ that are not in context $s$ for any $s\in\tau$:
\begin{align*}
I(z_{1:n}; \tau) &= \{i \in\{1,\dots,n\} : \forall s\in \tau, z_{(i-|s|)\vee 1:i-1}\neq s\} 
\;.\end{align*}
  Then, if $\P_{\tau,\theta}\left(Z_{1:n}=z_{1:n}\right)>0$,
\begin{multline*}
\label{model2}
\P_{\tau,\theta}\left(Z_{1:n}=z_{1:n}\right)=\prod_{i\in I(z_{1:n};\tau)}\P_{\tau,\theta}\left(Z_{i}=z_{i}\vert Z_{1:i-1}=z_{1:i-1}\right)\\
\prod_{s\in\tau}
P_{\theta(s)}\left(S(s;z_{1:n})\right) \;,
\end{multline*}
where for $\vartheta\in {\cal P}_{A}$,  $P_{\vartheta}$ denotes the probability distribution of the memoryless source on $A$ with parameter $\vartheta$.

Assume that $X=(X_{n})_{n\in\Z}$ and $Y=(Y_{n})_{n\in\Z}$ are independent stationary context tree sources. % with probability distributions $\P_{X}$ and $\P_{Y}$,  that share some contexts and associated conditional probability distributions, which means that
Let us define subsets $\sigma_{0}$, $\sigma_{1}$ and $\sigma_{2}$ of $A^{\star}$,
and parameters  $\theta_{0}=(\theta_{0}(s))_{s\in\sigma_{0}}$, $\theta_{1}=(\theta_{1}(s))_{s\in\sigma_{1}}$, $\theta_{2}=(\theta_{2}(s))_{s\in\sigma_{2}}$, $\theta_{i}(s)\in{\cal P}_{A}$, $s\in\sigma_{i}$, $i=0,1,2$ by the following properties: $X$ has distribution $\P_{\tau_{1},
( \theta_{0},\theta_{1})}$, $Y$ has distribution $\P_{\tau_{2}, (\theta_{0},\theta_{2})}$, and
\begin{equation}
\label{T1}
\sigma_{1}\cap \sigma_{0} = \emptyset,\;\sigma_{2}\cap \sigma_{0} = \emptyset,
\end{equation}
\begin{equation}
\label{T2}
\tau_{1}:=\sigma_{1}\cup \sigma_{0} \text{ is a complete tree},
\end{equation}
\begin{equation}
\label{T3}
\tau_{2}:=\sigma_{2}\cup \sigma_{0} \text{ is a complete tree},
\end{equation}
\begin{equation}
\label{T4}
\forall s\in \sigma_{1}\cap \sigma_{2},\;\theta_{1}\left(s\right)\neq \theta_{2}\left(s\right).
\end{equation}
$\sigma_{0}$ is the set of shared contexts, that is the set of contexts which intervene in both sources with the same associated probability distributions.

Given two samples  $X_{1:n}=(X_{1},\ldots,X_{n})$ and $Y_{1:m}=(Y_{1},\ldots,Y_{m})$ generated by $X$ and $Y$ respectively,
the aim of this paper is to propose a statistical method for the joint estimation of $\sigma_{0}$, $\sigma_{1}$ and $\sigma_{2}$,
and consequently of $\theta_{0}$, $\theta_{1}$ and $\theta_{2}$. 

This is a model selection problem, in which the collection of models is described by possible $\sigma_{0}$, $\sigma_{1}$ and $\sigma_{2}$,'s and for fixed $\sigma_{0}$, $\sigma_{1}$ and $\sigma_{2}$ the model consists of all 
$\P_{\sigma_{1}\cup \sigma_{0}, (\theta_{0},\theta_{1})}$ and $\P_{\sigma_{2}\cup \sigma_{0}, (\theta_{0},\theta_{2})}$ for 
any possible $\theta_{i}$, $i=0,1,2$.\\

We propose in the next section a selection method using penalized maximum likelihood for the entire set of observations.

%If one now gets  samples of independent sources $Z^{1},\ldots,Z^{N}$ that have $\P_{X}$ or $\P_{Y}$ distribution
%($\P_{X}$ and $\P_{Y}$ unknown but sharing contexts), and if one wants to find the shared contexts, and decide for each sample from which of both context source it comes from, one faces an unsupervised classification problem that can be modeled as a mixture. CAN WE DO SOMETHING ON THAT QUESTION ????

\section{The joint Context Tree Estimator}\label{sec:estimator}
\subsection{Likelihood in context-tree models}
For any $(\sigma_{0},\sigma_{1},\sigma_{2})$ satisfying (\ref{T1}), (\ref{T2}) and (\ref{T3}), define $\M_{(\sigma_{0},\sigma_{1},\sigma_{2})}$ as the set of distributions $\Q$ on $A^{\N}\times A^{\N}$ of form
$$
\Q=\P_{\sigma_{1}\cup \sigma_{0}, (\theta_{0},\theta_{1})}\otimes \P_{\sigma_{2}\cup \sigma_{0},( \theta_{0},\theta_{2})}:= \Q_{X}\otimes \Q_{Y}
$$
for some $\theta_{0}=(\theta_{0}(s))_{s\in\sigma_{0}}$, $\theta_{1}=(\theta_{1}(s))_{s\in\sigma_{1}}$, $\theta_{2}=(\theta_{2}(s))_{s\in\sigma_{2}}$, such that $\theta_{i}(s)\in{\cal P}_{A}$, $s\in\sigma_{i}$, $i=0,1,2$. Here we do not assume (\ref{T4}).\\
For any integers $n$ and $m$, any $x_{1:n}\in A^{n}$ and $y_{1:m}\in A^{m}$ and
any string $s$, denote by $S(s;x_{1:n};y_{1:m})=S(s;x_{1:n})S(s;y_{1:m})$ the concatenation of the $x_{i}$'s in context $s$, and of the
$y_{i}$'s in context $s$.
One has :
\begin{multline}
\label{model3}
\Q\left(X_{1:n}=x_{1:n}; Y_{1:m}=y_{1:m}\right)=\\
\prod_{i\in I(x_{1:n};\sigma_{1}\cup \sigma_{0})}\P_{\sigma_{1}\cup \sigma_{0}, (\theta_{0},\theta_{1})}\left(X_{i}=x_{i}\vert X_{1:i-1}=x_{1:i-1}\right)\\
\prod_{i\in I(y_{1:m};\sigma_{2}\cup \sigma_{0})}\P_{\sigma_{2}\cup \sigma_{0}, (\theta_{0},\theta_{2})}\left(Y_{i}=y_{i}\vert Y_{1:i-1}=y_{1:i-1}\right)
\\
\prod_{s\in\sigma_{0}}
P_{\theta_{0}(s)}\left(S(s;x_{1:n};y_{1:m})\right)\prod_{s\in\sigma_{1}}
P_{\theta_{1}(s)}\left(S(s;x_{1:n})\right)\prod_{s\in\sigma_{2}}
P_{\theta_{2}(s)}\left(S(s;y_{1:m})\right).
\end{multline}
Let us now note for any $s\in A^{\star}$ and any $a\in A$:
$$
N_{n,X}\left(s,a\right)=\sum_{i=|s|+1}^{n}\1_{X_{i-|s|:i-1}=s, X_{i}=a},\;
N_{n,X}\left(s\right)=\sum_{i=|s|+1}^{n}\1_{X_{i-|s|:i-1}=s}
$$
where it is understood that an empty sum is $0$,
and 
$$
N_{m,Y}\left(s,a\right)=\sum_{i=|s|+1}^{m}\1_{Y_{i-|s|:i-1}=s,Y_{i}=a},\;
N_{m,Y}\left(s\right)=\sum_{i=|s|+1}^{m}\1_{Y_{i-|s|:i-1}=s}.
$$
Observe that $N_{n,X}\left(\epsilon\right) = n$ and $N_{m,Y}\left(\epsilon\right) =m$.
Then, when maximizing over $\M_{(\sigma_{0},\sigma_{1},\sigma_{2})}$ the likelihood
as given by (\ref{model3}),
we shall use the approximation that the first two terms may be maximized as free parameters (so that their maximization gives $1$). Thus we shall use the pseudo maximum log-likelihood
\begin{align*}
\ell_{n,m}\Big(&\sigma_{0},\sigma_{1},\sigma_{2}\Big)= \sum_{s\in\sigma_{1}}\sum_{a\in A} N_{n,X}\left(s,a\right)\log\left(\frac{N_{n,X}\left(s,a\right)}{N_{n,X}\left(s\right)}\right)\\&+
\sum_{s\in\sigma_{2}}\sum_{a\in A} N_{m,Y}\left(s,a\right)\log\left(\frac{N_{m,Y}\left(s,a\right)}{N_{m,Y}\left(s\right)}\right)
\\
&+
\sum_{s\in\sigma_{0}}\sum_{a\in A} \left[N_{n,X}\left(s,a\right)+N_{m,Y}\left(s,a\right)\right]\log\left(\frac{N_{n,X}\left(s,a\right)+N_{m,Y}\left(s,a\right)}{N_{n,X}\left(s\right)+N_{m,Y}\left(s\right)}\right),
\end{align*}
where by convention for any non negative integer $p$, $0\log \frac{0}{p}=0$. Here $\log u$ denotes the logarithm of $u$ in base $2$.\\
For any string $s$, we shall write $Q_{X}\left(\cdot \vert s\right)$ and $Q_{Y}\left(\cdot \vert s\right)$ the probability distributions on $A$ given by: $\forall a\in A,$
\begin{align*}
 &Q_{X}\left(a \vert s\right)=\Q\left(X_{|s|+1}=a \vert X_{1:|s|}=s\right),\;\\
&Q_{Y}\left(a \vert s\right)=\Q\left(Y_{|s|+1}=a \vert Y_{1:|s|}=s\right),
\end{align*}
and $\widehat{Q}_{X}\left(\cdot \vert s\right)$, $\widehat{Q}_{Y}\left(\cdot \vert s\right)$ and $\widehat{Q}_{XY}\left(\cdot \vert s\right)$ the probability distributions on $A$ given by: $\forall a\in A$
\begin{align*}
&\widehat{Q}_{X}\left(a \vert s\right)=\frac{N_{n,X}\left(s,a\right)}{N_{n,X}\left(s\right)}\;,\qquad \widehat{Q}_{Y}\left(a \vert s\right)=\frac{N_{m,Y}\left(s,a\right)}{N_{m,Y}\left(s\right)}\\
&\widehat{Q}_{XY}\left(a \vert s\right)=\frac{N_{n,X}\left(s,a\right)+N_{m,Y}\left(s,a\right)}{N_{n,X}\left(s\right)+N_{m,Y}\left(s\right)}
\end{align*}
whenever $N_{n,X}(s)>0$, $N_{m, Y}(s)>0$ and $N_{n,X}(s)+N_{m, Y}(s)>0$ respectively.
In the same way, with some abuse of notation, we note $Q_{X}$ and $Q_{Y}$ any $|s|$-marginal probability distributions on $A^{|s|}$ defined respectively by $\Q_{X}$ and $\Q_{Y}$.
\\
\subsection{Definition of the joint estimator}
Let $\pen(\cdot)$ be a function from $\N$ to $\R$, % possibly also of the observations $X_{1:n}$ and $Y_{1:m}$, 
which will be called penalty function, and define the estimators $\widehat{\sigma}_{0}$, $\widehat{\sigma}_{1}$ and $\widehat{\sigma}_{2}$ as a triple of maximizers of
\begin{multline*}
C_{n,m}\left(\sigma_{0},\sigma_{1},\sigma_{2}\right)=\ell_{n,m}\left(\sigma_{0},\sigma_{1},\sigma_{2}\right)\\-\frac{(|A|-1)}{2}(|\sigma_{0}|\pen(n+m)+|\sigma_{1}|\pen(n)+|\sigma_{2}|\pen(m))
\end{multline*}
over all possible $(\sigma_{0},\sigma_{1},\sigma_{2})$ satisfying (\ref{T1}), (\ref{T2}) and (\ref{T3}). The BIC estimator corresponds to the choice $\pen(\cdot) = \log(\cdot)$.
Notice that it is enough to restrict the maximum over sets $\sigma_{0},\sigma_{1},\sigma_{2}$ that have strings $s$ with length $|s|\leq n\vee m -1$.
Indeed, if a string $s$ has length $|s|\geq n$, then for any $a\in A$, $N_{n,X}(s,a)=0$,
if  $s$ has length $|s|\geq m$, then for any $a\in A$, $N_{m,Y}(s,a)=0$.
\\
For any integer $D$, denote  
$$
\left(\widehat{\sigma}_{D,0}, \widehat{\sigma}_{D,1},\widehat{\sigma}_{D,2}\right)=\argmax
C_{n,m}\left(\sigma_{0},\sigma_{1},\sigma_{2}\right)
$$
where the maximization is over all $(\sigma_{0},\sigma_{1},\sigma_{2})$ satisfying (\ref{T1}), (\ref{T2}) and (\ref{T3}) and such that
for any $s\in \sigma_{0}\cup \sigma_{1}\cup \sigma_{2}$, $|s|\leq D$. 
Then, as explained before, the joint estimator $\left(\widehat{\sigma}_{0}, \widehat{\sigma}_{1},\widehat{\sigma}_{2}\right)$ is seen to be:
$$
\left(\widehat{\sigma}_{0}, \widehat{\sigma}_{1},\widehat{\sigma}_{2}\right)=
\left(\widehat{\sigma}_{n\vee m -1,0}, \widehat{\sigma}_{n\vee m -1,1},\widehat{\sigma}_{n\vee m -1,2}\right)\;.
$$

\subsection{Consistency of the joint estimator}
Now assume that $X$ and $Y$ are independent with distribution
$$
\Q^{\star}=\P_{\sigma_{1}^{\star}\cup \sigma_{0}^{\star}, (\theta_{0}^{\star},\theta_{1}^{\star})}\otimes \P_{\sigma_{2}^{\star}\cup \sigma_{0}^{\star}, (\theta_{0}^{\star},\theta_{2}^{\star})}
$$
where $\sigma_{0}^{\star}$, $\sigma_{1}^{\star}$, $\sigma_{2}^{\star}$ are finite subsets of $A^{\star}$
satisfying (\ref{T1}), (\ref{T2}) and (\ref{T3}), and such that (\ref{T4}) holds. 
%Let $D^{\star}$ be the maximum length of strings in $\sigma_{0}^{\star}\cup \sigma_{1}^{\star}\cup \sigma_{2}^{\star}$:
%$$
%D^{\star}=\max \left\{|s|\;:\;s\in \sigma_{0}^{\star}\cup \sigma_{1}^{\star}\cup \sigma_{2}^{\star} \right\}.
%$$
%Define for any integer $D$:
%$$
%\left(\widehat{\sigma}_{D,0}, \widehat{\sigma}_{D,1},\widehat{\sigma}_{D,2}\right)=\argmax
%C_{n,m}\left(\sigma_{0},\sigma_{1},\sigma_{2}\right)
%$$
%over all $(\sigma_{0},\sigma_{1},\sigma_{2})$ satisfying (\ref{T1}), (\ref{T2}) and (\ref{T3}) and such that
%for any $s\in \sigma_{0}\cup \sigma_{1}\cup \sigma_{2}$, $|s|\leq D$.
\begin{theo}
\label{th:consi}
Assume that $n$ and $m$ go to infinity in such a way that
\begin{equation}
\label{equilibre}
\lim_{n\rightarrow \infty}\frac{n}{m}=c,\;0<c<+\infty.
\end{equation}
Assume 
moreover that for any integer $k$,
%the penalty function verifies 
%$$
%\lim_{n\rightarrow \infty}\frac{\pen\left(n\right)}{n}=0,\;
%\limsup_{n\rightarrow \infty}\frac{\log n}{\pen\left(n\right)}< +\infty.
%$$
$$
\pen\left(k\right)=\log k.
$$
Then the joint estimator is consistent, i.e.
$$
\left(\widehat{\sigma}_{0}, \widehat{\sigma}_{1},\widehat{\sigma}_{2}\right)=
\left(\sigma_{0}^{\star}, \sigma_{1}^{\star}, \sigma_{2}^{\star} \right)
$$
$\Q^{\star}$-eventually almost surely as $n$ goes to infinity.
\end{theo}
We have presented our joint estimator with a generic penalty $\pen (\cdot)$, and Section \ref{sec:algo} describes a procedure for computing efficiently this estimator in the general case.
However, the  consistency result only covers the choice of the BIC penalty \cite{schwarz1978}, that is the penalty which is the logarithm of the number of observations times half the number of free parameters.
The proof of Theorem~\ref{th:consi} is given in Section~\ref{sec:proof}.

\section{An Efficient algorithm for the joint estimator}\label{sec:algo}
In this section, we propose an efficient algorithm for the computation of the joint estimator with no restriction on the depth of the trees.
The recursive tree structure makes it possible to maximize the penalized maximum likelihood criterion without considering all possible models (which are far too numerous).
The greedy algorithm we present here can be seen as a non-trivial extension of the Context Tree Maximization algorithm that was first presented in \cite{willems1995}, see also \cite{csiszar2006}.
For each possible node $s$ of the estimated tree, the algorithm first computes recursively, from the leaves to the root, indices $\chi_{s} (X_{1:n}),\chi_{s} (Y_{1:m})$ and $\chi_{s} (X_{1:n};Y_{1:m})$. In a second step, the estimated tree is constructed from the root to the leaves according to these indices.

For any string $s$ let
\begin{align*}
&\widehat{P}_{s}\left(X_{1:n}\right) = \prod_{a\in A}\left(\frac{N_{n,X}\left(s,a\right)}{N_{n,X}\left(s\right)}\right)^{N_{n,X}\left(s,a\right)} \;,\\
& \widehat{P}_{s}\left(Y_{1:m}\right) = \prod_{a\in A}\left(\frac{N_{m,Y}\left(s,a\right)}{N_{m,Y}\left(s\right)}\right)^{N_{m,Y}\left(s,a\right)},
\end{align*}

and let
$$
\widehat{P}_{s}\left(X_{1:n};Y_{1:m}\right) = \prod_{a\in A}\left(\frac{N_{n,X}\left(s,a\right)+N_{m,Y}\left(s,a\right)}{N_{n,X}\left(s\right)+N_{m,Y}\left(s\right)}\right)^{N_{n,X}\left(s,a\right)+N_{m,Y}\left(s,a\right)}
$$
where again it is understood that for any non negative integer $n$, $(\frac{0}{n})^{0}=1$. Notice that, because of possible side effects,  $\widehat{P}_{s}\left(X_{1:n};Y_{1:m}\right)$ is not in general equal to $\widehat{P}_{s}\left(X_{1:n}Y_{1:m}\right)$.\\

\subsubsection*{Step 1: computation of the indices}
For any set of strings $\sigma$, we denote by $\sigma s$  the set of strings $us$, $u\in \sigma$: $\sigma s= \{us: u\in\sigma\}$. 
Let $\sigma$ be a tree, and let
\def\cout{R}
$$
\cout_{\sigma ; s}\left(X_{1:n}\right) = \sum_{u\in\sigma s} \log \widehat{P}_{u}\left(X_{1:n}\right)- |\sigma| \pen \left(n\right),
$$
$$
\cout_{\sigma ; s}\left(Y_{1:m}\right) = \sum_{u\in\sigma s} \log \widehat{P}_{u}\left(Y_{1:m}\right)- |\sigma| \pen \left(m\right),
$$
$$
\cout_{\sigma ; s}\left(X_{1:n};Y_{1:m}\right) = \sum_{u\in\sigma s} \log \widehat{P}_{u}\left(X_{1:n};Y_{1:m}\right)- |\sigma| \pen \left(n+m\right).
$$
Let $D$ be an upper-bound on the size of the candidate contexts in $\sigma_0\cup\sigma_1\cup\sigma_2$. 
Note that it is sufficient to consider $D=n\vee m$ to investigate all possible trees.
Define for any string of length $|s|=D$:
$$
V_{s}\left(X_{1:n}\right)=\cout_{\{\epsilon\} ; s}\left(X_{1:n}\right),\;\chi_{s}\left(X_{1:n}\right)=0,$$
$$
V_{s}\left(Y_{1:m}\right)=\cout_{\{\epsilon\} ; s}\left(Y_{1:m}\right),\;\chi_{s}\left(Y_{1:m}\right)=0,$$
$$
V_{s}\left(X_{1:n};Y_{1:m}\right)=\max\left\{\cout_{\{\epsilon\};s}\left(X_{1:n};Y_{1:m}\right);\cout_{\{\epsilon\} ; s}\left(X_{1:n}\right)+\cout_{\{\epsilon\}; s}\left(Y_{1:m}\right)
\right\},
$$
and 
$$\chi_{s}\left(X_{1:n};Y_{1:m}\right)=\left\{
\begin{array}{l}
1\;{\text{, if }} \;V_{s}\left(X_{1:n};Y_{1:m}\right)=\cout_{\{\epsilon\};s}\left(X_{1:n};Y_{1:m}\right)\\
2\;{\text{, else}}.
\end{array}
\right.
$$
Then compute recursively for all $s$ such that $|s|<D$:
$$
V_{s}\left(X_{1:n}\right)=\max\left\{\cout_{\{\epsilon\} ; s}\left(X_{1:n}\right);
\sum_{a\in A}V_{as}\left(X_{1:n}\right)\right\},
$$
and
$$\chi_{s}\left(X_{1:n}\right)=\left\{
\begin{array}{l}
0\;{\text{, if }} \;V_{s}\left(X_{1:n}\right)=\cout_{\{\epsilon\};s}\left(X_{1:n}\right)\\
1\;{\text{ else}},
\end{array}
\right.
$$
$$
V_{s}\left(Y_{1:m}\right)=\max\left\{\cout_{\{\epsilon\} ; s}\left(Y_{1:m}\right);
\sum_{a\in A}V_{as}\left(Y_{1:m}\right)\right\},
$$
and
$$\chi_{s}\left(Y_{1:m}\right)=\left\{
\begin{array}{l}
0\;{\text{, if }} \;V_{s}\left(Y_{1:m}\right)=\cout_{\{\epsilon\};s}\left(Y_{1:m}\right)\\
1\;{\text{ else}}.
\end{array}
\right.
$$
Define also
$$
V_{s}\left(X_{1:n};Y_{1:m}\right)=\max\left\{
\begin{array}{l}
\cout_{\{\epsilon\};s}\left(X_{1:n};Y_{1:m}\right)\\
V_{s}\left(X_{1:n}\right)+V_{s}\left(Y_{1:m}\right)\\
\sum_{a\in A}V_{as}\left(X_{1:n};Y_{1:m}\right),
\end{array}\right.
$$
and 
$$
\chi_{s}\left(X_{1:n};Y_{1:m}\right)=\left\{
\begin{array}{l}
1\;{\text{, if }}\;V_{s}\left(X_{1:n};Y_{1:m}\right)= \cout_{\{\epsilon\};s}\left(X_{1:n};Y_{1:m}\right)\;,\\
2\;{\text{, if }}\;V_{s}\left(X_{1:n};Y_{1:m}\right)= 
V_{s}\left(X_{1:n}\right)+V_{s}\left(Y_{1:m}\right)\;,\\
3\;{\text{ else. }}
\end{array}\right.
$$
For any $(\sigma_{0},\sigma_{1},\sigma_{2})$ satisfying (\ref{T1}), (\ref{T2}) and (\ref{T3}), define
$$
\cout_{(\sigma_{1},\sigma_{2},\sigma_{0}) ; s}\left(X_{1:n};Y_{1:m}\right) = \cout_{\sigma_{1} ; s}\left(X_{1:n}\right)+\cout_{\sigma_{2} ; s}\left(Y_{1:m}\right) +
\cout_{\sigma_{0} ; s}\left(X_{1:n};Y_{1:m}\right).
$$
Notice that
$$
\cout_{(\sigma_{1},\sigma_{2},\emptyset) ; s}\left(X_{1:n};Y_{1:m}\right)=
\cout_{\sigma_{1}; s}\left(X_{1:n}\right)+\cout_{\sigma_{2}; s}\left(Y_{1:m}\right)
$$
and
$$
\cout_{(\emptyset,\emptyset,\sigma_{0}) ; s}\left(X_{1:n};Y_{1:m}\right)=\cout_{\sigma_{0} ; s}\left(X_{1:n};Y_{1:m}\right).
$$
Moreover, remark that
\begin{itemize}
 \item either $\sigma_{1}$ and $\sigma_{2}$ are the empty set and $\sigma_{0}$ is not the empty set,
 \item or $\sigma_{0}$ is the empty set and neither $\sigma_1$ nor $\sigma_2$ are the empty set,
 \item or none of them is the empty set.
\end{itemize}

\subsubsection*{Step 2: construction of the estimated trees}
Once the indicators $\chi_{s}\left(X_{1:n}\right)$ and $\chi_{s}\left(X_{1:n}\right)$ have been computed, the estimated sets can be computed recursively from the root to the leaves.
Recall that Csiszar and Talata \cite{csiszar2006} prove that for any string $s$ such that $|s|\leq D$:
\begin{equation}
V_{s}\left(X\right)=\max_{\sigma} \cout_{\sigma;s}\left(X\right)
\label{eq:vx}\end{equation}
and 
\begin{equation}
V_{s}\left(Y\right)=\max_{\sigma} \cout_{\sigma;s}\left(Y\right).
\label{eq:vy}\end{equation}
Call $\sigma_{X_{1:n}}\left(s \right)$ (resp. $\sigma_{Y_{1:m}}\left(s \right)$) a tree maximizing \eqref{eq:vx} (resp. \eqref{eq:vy}).
$\sigma_{X_{1:n}}\left(s \right)$ and $\sigma_{Y_{1:m}}\left(s \right)$ can be computed recursively as follows:
start with the strings $s$ of length $D$;
\begin{itemize}
\item if $\chi_{s}\left(X_{1:n}\right)=0$, then $\sigma_{X_{1:n}}\left(s \right)=\{\epsilon\}$,
\item if $\chi_{s}\left(X_{1:n}\right)=1$, then $\sigma_{X_{1:n}}\left(s \right)=\cup_{a\in A}\sigma_{X_{1:n}}\left(as \right) a$,
\item if $\chi_{s}\left(Y_{1:m}\right)=0$, then $\sigma_{Y_{1:m}}\left(s \right)=\{\epsilon\}$,
\item
if $\chi_{s}\left(Y_{1:m}\right)=1$, then $\sigma_{Y_{1:m}}\left(s \right)=\cup_{a\in A}\sigma_{Y_{1:m}}\left(as \right) a$.
\end{itemize}

Namely, for any string $s$ such that $|s|\leq D$, define  $\sigma_{1}\left(s \right)$, 
$\sigma_{2}\left(s \right)$ and $\sigma_{0}\left(s \right)$ as: 
\begin{itemize}
\item
if $\chi_{s}\left(X_{1:n};Y_{1:m}\right)=1$, then $\sigma_{1}\left(s \right)=\sigma_{2}\left(s \right)=\emptyset$
and $\sigma_{0}\left(s \right)=\{\epsilon\}$,
\item
if $\chi_{s}\left(X_{1:n};Y_{1:m}\right)=2$, then $\sigma_{1}\left(s \right)=\sigma_{X_{1:n}}\left(s \right)$, $\sigma_{2}\left(s \right)=\sigma_{Y_{1:m}}\left(s \right)$
and $\sigma_{0}\left(s \right)=\emptyset$,
\item
if $\chi_{s}\left(X_{1:n};Y_{1:m}\right)=3$, then $\sigma_{1}\left(s \right)=\cup_{a\in A}\sigma_{1}\left(as \right) a$,
$\sigma_{2}\left(s \right)=\cup_{a\in A}\sigma_{2}\left(as \right) a$
and $\sigma_{0}\left(s \right)=\cup_{a\in A}\sigma_{0}\left(as \right) a$.
\end{itemize}
\subsubsection*{Validity of the algorithm}
The next proposition shows that the two-step procedure described above computes the maximum pseudo-likelihood estimator in the joint model.
\begin{prop}
\label{computree}
For any string $s$ such that $|s|\leq D$,
$$
V_{s}\left(X_{1:n};Y_{1:m}\right)=\max \cout_{(\sigma_{1},\sigma_{2},\sigma_{0}) ; s}\left(X_{1:n};Y_{1:m}\right)
$$
where the maximum is over all $(\sigma_{0},\sigma_{1},\sigma_{2})$ that verify (\ref{T1}), (\ref{T2}) and (\ref{T3})
and such that
$$\forall u \in \sigma_{1}\cup \sigma_{2} \cup \sigma_{0},\; |u| + |s| = D.$$
In particular,
$$\widehat{\sigma}_{D,0}=\sigma_{0}\left(\epsilon \right),\;\widehat{\sigma}_{D,1}=\sigma_{1}\left(\epsilon \right), \;\widehat{\sigma}_{D,2}=\sigma_{2}\left(\epsilon \right).$$
\end{prop}
{\bf{Proof:}}\\
The proof is by induction. Observe first that
$$
V_{s}\left(X_{1:n}\right)+V_{s}\left(Y_{1:m}\right)=\max_{\sigma_{1},\sigma_{2}} \cout_{(\sigma_{1},\sigma_{2},\emptyset) ; s}\left(X_{1:n};Y_{1:m}\right).
$$
Now,
if $|s|=D$, then either $\sigma_{1}=\sigma_{2}=\{\epsilon\}$ and $\sigma_{0} = \emptyset$, or
$\sigma_{1}=\sigma_{2}=\emptyset$ and $\sigma_{0} =\{\epsilon\}$, and we have
$$
V_{s}\left(X_{1:n};Y_{1:m}\right)=\max \left\{\cout_{(\{\epsilon\},\{\epsilon\},\emptyset) ; s}\left(X_{1:n};Y_{1:m}\right);
\cout_{(\emptyset,\emptyset,\{\epsilon\}) ; s}\left(X_{1:n};Y_{1:m}\right)
\right\}.
$$
Let us now take $|s|<D$ and assume that Proposition \ref{computree} is true for all strings $as$, $a\in A$. 
The maximum of the  $\cout_{(\sigma_{1},\sigma_{2},\sigma_{0}) ; s}\left(X_{1:n};Y_{1:m}\right)$ over all $(\sigma_{0},\sigma_{1},\sigma_{2})$ that verify (\ref{T1}), (\ref{T2}) and (\ref{T3}) and such that $\forall u \in \sigma_{1}\cup \sigma_{2} \cup \sigma_{0},\; |u| + |s| = D$, is reached by a triple $(\sigma_{1}, \sigma_{2}, \sigma_{0})$ such that:
\begin{itemize}
\item either $\sigma_0=\{\epsilon\}$, in which case $\sigma_1$ and $\sigma_2$ are necessarily empty and \[  \cout_{(\sigma_{1},\sigma_{2},\sigma_{0});s}\left(X_{1:n};Y_{1:m}\right)  = \cout_{(\emptyset,\emptyset,\{\epsilon\});s}\left(X_{1:n};Y_{1:m}\right) = \cout_{\{\epsilon\};s}\left(X_{1:n};Y_{1:m}\right);\] 
\item or at least one among $\sigma_1$ and $\sigma_2$ is equal to $\{\epsilon\}$: then $\sigma_{0}=\emptyset$ and 
\[ \cout_{(\sigma_{1},\sigma_{2},\sigma_{0});s}\left(X_{1:n};Y_{1:m}\right) =  \cout_{\sigma_1;s}(X_{1:n}) + \cout_{\sigma_2;s}(Y_{1:m}) = V_s(X_{1:n}) + V_s(Y_{1:m})\]
as in \cite{csiszar2006};
\item or $\sigma_{1}, \sigma_{2}, \sigma_{0}$ are all different from $\{\epsilon\}$, and then each $\sigma_i, 0\leq i\leq 2$ can be written as 
$\sigma_i= \cup_{a\in A}\sigma_i(a) a$;  note that it is possible that, for some $i\in\{0,1,2\}$ and some $a\in A$, $\sigma_i(a)$ is empty, or even that $\sigma_i$ is empty. In any case, for each $a\in A$ it is easily checked that  $\sigma_1(a), \sigma_2(a)$ and $\sigma_0(a)$ satisfy (\ref{T1}),  (\ref{T2}) and (\ref{T3}). Thus
\begin{align*}
\cout_{(\sigma_{1},\sigma_{2}, \sigma_{0});s}\left(X_{1:n};Y_{1:m}\right) &= \sum_{a\in A}
\cout_{(\sigma_{1}(a),\sigma_{2}(a), \sigma_{0}(a));as}\left(X_{1:n};Y_{1:m}\right)\\
&= \sum_{a\in A}
\max_{\bar{\sigma}_{1},\bar{\sigma}_{2},\bar{\sigma}_{0}}
\cout_{(\sigma_{1},\sigma_{2}, \sigma_{0});as}\left(X_{1:n};Y_{1:m}\right)\\
& = \sum_{a\in A}V_{as}\left(X_{1:n};Y_{1:m}\right)
\end{align*}
by the induction hypothese.
\end{itemize}
To conclude the proof, it is enough to be reminded that, by definition, 
\begin{multline*}
V_{s}\left(X_{1:n};Y_{1:m}\right) = \max\Bigg\{
\cout_{\{\epsilon\};s}\left(X_{1:n};Y_{1:m}\right),\\
V_{s}\left(X_{1:n}\right)+V_{s}\left(Y_{1:m}\right), \quad
\sum_{a\in A}V_{as}\left(X_{1:n};Y_{1:m}\right)\Bigg\}\;.
\end{multline*}

Obviously the computational complexity of this procedure is proportional to the number of candidate nodes $s$, which is equal to the number of distinct subsequences of $X_{1:n}$ and $Y_{1:m}$, and hence quadratic in $n$ and $m$. However, if necessary, it is possible to obtain a linear complexity algorithm by using compact suffix trees, as explained in \cite{garivier2006}.

% \item 
%  \item either $\sigma_{1}=\sigma_{2}=\emptyset$, $\sigma_{0} \neq \emptyset$ and $\sigma_{0} \neq \{\epsilon\}$,
%  \item or $\sigma_{0} = \emptyset$ and $\sigma_{1}$, $\sigma_{2}$ are neither empty nor equal to $\{\epsilon\}$
%  \item or all trees neither empty and nor equal to $\{\epsilon\}$.
% \end{itemize}
% Thus, the definition of $V_{s}\left(X_{1:n};Y_{1:m}\right)$ and the induction hypothese imply that
% \begin{align*}
% V_{s}\left(X_{1:n};Y_{1:m}\right) &=\max\left\{
% \begin{array}{l}
% \cout_{(\emptyset,\emptyset,\{\epsilon\});s}\left(X_{1:n};Y_{1:m}\right)\\
% \max_{\sigma_{1},\sigma_{2}} \cout_{(\sigma_{1},\sigma_{2},\emptyset) ; s}\left(X_{1:n};Y_{1:m}\right)\\
% \max_{\sigma_{1}\neq \{\epsilon\},\sigma_{2}\neq \{\epsilon\},\sigma_{0}\neq \{\epsilon\}}
% \cout_{(\sigma_{1},\sigma_{2}, \sigma_{0});s}\left(X_{1:n};Y_{1:m}\right)
% \end{array}\right.\\
%  &=\max \cout_{(\sigma_{1},\sigma_{2},\sigma_{0}) ; s}\left(X_{1:n};Y_{1:m}\right).
%\end{align*}

\section{Simulation study}
\label{sec:simu}
In this section, we experimentally show the value of joint estimation when the two sources $X$ and
$Y$ share some contexts. 
We compare the results obtained by the BIC joint-estimator described above with the following direct approach.
First, we estimate $\tau_X$ using the standard BIC tree estimate $\hat{\tau}_X= \hat{\tau}_X(X_{1:n})$, and we independently estimate $\tau_Y$ using $\hat{\tau}_Y= \hat{\tau}_Y(Y_{1:m})$.
Then, for all contexts $s$ that are present both in $\hat{\tau}_X$ and in $\hat{\tau}_Y$,  we compute the chi-squared distance of the conditional empirical distributions: if this distance is smaller than a given threshold, we decide that $s$ is a shared context. The value of the threshold was chosen in order to maximize the frequency of correct estimation.

\subsection{A particularly favorable example}

First consider the following case:
\begin{itemize}
\item $X$ and $Y$ are $\{1,2\}$-valued context-tree sources;
%\item $\P_X$ is defined by the conditional distributions $P(1|1) = 1/3, P(1|12)=1/3, P(1|22)=2/3$;
%\item $\P_Y$ is defined by the conditional distributions $P(1|1) = 1/3, P(1|12)=1/3, P(1|122)=1/5, P(1|222)=1/2$;
%\item the estimates are computed from $X_{1:n}$ and $Y_{1:m}$ with $n=150$ and $m=300$;
\item $\Q_X$ is defined by the conditional distributions $\Q_X(X_0=1|X_{-1}=1) = 1/3, \Q_X(X_0=1|X_{-2:-1}=12)=1/3, \Q_X(X_0=1|X_{-2:-1}=22)=2/3$;
\item $\Q_Y$ is defined by the conditional distributions $\Q_Y(Y_0=1|Y_{-1}=1) = 3/4, \Q_Y(Y=0=1|Y_{-2:-1}=12)=1/3, \Q_Y(Y_0=1|Y_{-2:-1}=22)=2/3$;
\item the estimates are computed from $X_{1:n}$ and $Y_{1:m}$ with $n=500$ and $m=1000$;
\item the probability of correctly identifying the tree by each method is estimated by a Monte-Carlo procedure with $1000$ replications (margin of error $\approx 1.5\%$).
\end{itemize}
%In the next section, we describe an efficient algorithm that permits to compute the joint context tree estimator without computing the penalized criterion for all possible models. 
%
%It appears that separate estimates for $X$ (resp. $Y$) are correct with probability about $65\%$ (resp. $38\%$, while
%using joint estimation the trees are correctly identified with probability respectively $53\%$ and $67\%$. 
%More interestingly, the probability that \emph{both} estimated trees are correct raises from $26\%$ to $40\%$.
In that example, we hence have $\sigma_0=\{12, 22\}$,  $\sigma_1=\{1\}$ and $\sigma_2 = \{1\}$.
We compare our joint estimation procedure with separate estimation using the following criteria: 
\begin{itemize}
\item the probability of correctly identifying $\tau_X$ (resp. $\tau_Y$);
\item the probability of correctly identifying simultaneously $\tau_X$ and $\tau_Y$;
\item the probability of correctly identifying $\sigma_0,\sigma_1,\sigma_2$;
\item the Kullback-Leibler divergence rates $\kl(\Q_Z|\hat{\Q}_Z)$  between the stationary processes $\Q_Z$ and $\hat{\Q}_Z$ for $Z\in\{X,Y\}$, which are computated  by using the fact that both $X$ and $Y$ are Markov chains of finite order.
\end{itemize}
The results are summarized in Figure~\ref{fig:resultsSimu1}.
It appears that the joint estimation approach has a significant advantage over separate estimation on all the criteria considered here, with one restriction: in some cases, the estimation of either $\tau_X$ or $\tau_Y$ can be deteriorated, while the other is (more significantly) improved.
In all cases, the probability of correctly estimating both $\tau_X$ and $\tau_Y$ at the same time is increased.
\begin{figure}
\begin{tabular}{c|cccccccc}
 & $\tau_X$ & $\tau_Y$ & $\tau_X$ and $\tau_Y$ & $ \sigma_0$ & $ \sigma_1$ & $ \sigma_2$  & $KL_X$ & $KL_Y$\\ \hline
sep. est. & $51\%$ & $44\%$ & $22\%$ &   $20\%$  &  $31\%$ &  $31\%$ &  $6.7\;10^{-3}$  &  $5.7\; 10^{-3}$   \\
joint est. &  $80\%$  &  $78\%$ & $76\%$ & $77\%$  &  $90\%$ &  $90\%$ &  $3.2\;10^{-3}$  &  $2.3\; 10^{-3}$ \\ 
\end{tabular}
\caption{Comparative performance of separate and joint estimation in a favorable case (probabilities of correct estimation). 
$KL_X$  and $KL_Y$ denote $\kl(\Q_X| \hat{\Q}_X)$ and $\kl(\Q_Y| \hat{\Q}_Y)$, respectively.}
\label{fig:resultsSimu1}
\end{figure}

\subsection{A less favorable example}
On the other hand, when $X$ and $Y$ share no (or few) contexts, then the joint estimation procedure can obviously only deteriorate the separate estimates by introducing some confusion between similar, but distinct conditional distributions of $X$ and $Y$.
An example of such a case is the following:  

\begin{itemize}
\item $X$ and $Y$ are $\{1,2\}$-valued context-tree sources;
%\item $\P_X$ is defined by the conditional distributions $P(1|1) = 1/3, P(1|12)=1/3, P(1|22)=2/3$;
%\item $\P_Y$ is defined by the conditional distributions $P(1|1) = 1/3, P(1|12)=1/3, P(1|122)=1/5, P(1|222)=1/2$;
%\item the estimates are computed from $X_{1:n}$ and $Y_{1:m}$ with $n=150$ and $m=300$;
\item $\Q_X$ is defined by the conditional distributions $\Q_X(X_0=1|X_{-1}=1) = 1/2, \Q_X(X_0=1|X_{-1}=2)=2/3$;
\item $\Q_Y$ is defined by the conditional distributions $\Q_Y(Y_0=1|Y_{-1}=1) = 1/2, \Q_Y(Y_0=1|Y_{-2:-1}=12)=3/5, \Q_Y(Y_0=1|Y_{-2:-1}=22)=3/4$;
\item the estimates are computed from $X_{1:n}$ and $Y_{1:m}$ with $n=1000$ and $m=1500$;
\item the probability of correctly identifying the tree by each method is estimated by a Monte-Carlo procedure with $1000$ replications (margin of error $\approx 1.5\%$).
\end{itemize}
In that example, $\sigma_0=\{1\}$,  $\sigma_1=\{2\}$ and $\sigma_2 = \{12,22\}$.
The results are summarized in Figure~\ref{fig:resultsSimu2}. In this case, $\Q_X$ and $\Q_Y$ are quite close, and the joint estimation procedure tends to merge them into a single, common distribution. Thus, the probability of correctly inferring the structure of $\Q_X$ and $\Q_Y$ is significantly deteriorated. 
%However, on may observe that the  Kullback-Leibler divergence between the stationary processes $\P_X$ and $\hat{P}_X$ (resp. $\P_Y$ and $\hat{P}_Y$) is less significant: the estimates are 
\begin{figure}
\begin{tabular}{c|cccccccc}
& $\tau_X$ & $\tau_Y$ & $\tau_X$ and $\tau_Y$ & $ \sigma_0$ & $ \sigma_1$ & $ \sigma_2$  & $KL_X$ & $KL_Y$\\ \hline
sep. est. & $97\%$ & $89\%$ & $86\%$ &   $84\%$  &  $84\%$ &  $82\%$ &  $1.0\;10^{-3}$  &  $1.3\; 10^{-3}$   \\
joint est. &  $60\%$  &  $76\%$ & $39\%$ & $40\%$  &  $40\%$ &  $39\%$ &  $1.7\;10^{-3}$  &  $2.0\; 10^{-3}$ \\ 
\end{tabular}

\caption{Comparative performance of separate and joint estimation in the unfavourable case (probabilities of correct estimation). 
$KL_X$  and $KL_Y$ denote $\kl(\Q_X| \hat{\Q}_X)$ and $\kl(\Q_Y| \hat{\Q}_Y)$, respectively..}
\label{fig:resultsSimu2}
\end{figure}

\subsection{Influence of the penalty term}

A natural question is whether the performance of joint (or even separate) estimation can be significantly improved by using other choices of penalty functions, especially choices of the form $\pen(n) = \lambda \log(n)$ for some positive $\lambda$. The BIC choice $\lambda=1$ may be improved by using a recent data-driven procedure called \emph{slope heuristic}, see \cite{birgeMassart07minipen}. However, in the present case, the attempts to tune the penalty function by using the slope heuristic merely resulted in a confirmation that the BIC choice could not be significantly improved on the examples considered here. In fact, in addition to the difficulty to detect the dimension gap and thus the minimal penalty in our simulations (which could be expected, as the number of models is very large whereas the sample are not huge), the ideal penalty estimator was never observed to be very different from $\lambda = 1$.

\subsection{Discussion}
The simulation study strongly indicates that the joint estimation procedure has a significantly improved performance when the two sources do share contexts and conditional distributions which appear with a significant probability in the samples.
On the other hand, when the sources share no or few contexts, the procedure may introduce some confusion between the estimates, as could be expected. 
%Better if there are shared contexts, but can possibly introduce confusion.

When the goal is joint estimation, deterioration in the estimation of one of the trees seems to be the price to pay for better estimating the other tree, and the net effect is positive. 

The predictive power of the estimated model is reflected by a measure of discrepancy between the true law of the process and the law of the estimated distribution. We chose to consider Kullback-Leibler divergence, as it is naturally associated to logarithmic prediction loss in information theory. As expected, a significant improvement is observed for the joint estimator in presence of shared contexts.
 
% for nx=200 and ny=400, we find:
%               TX -  TY -  TX and TY correctly estimated
%    separate : 80% - 55% - 44%
%    joint :    64% - 80 % - 55 %

% 
% \section{Comments and concluding remarks}
% \label{secconclu}
% 
% If one knows an upper bound on the depth of the context trees, then one may use a weaker penalty function, as in \cite{finesso1996}. More precisely,
% Assume $D \geq D^{\star}$, and that the penalty function verifies 
% $$
% \lim_{n\rightarrow \infty}\frac{\pen\left(n\right)}{n}=0,\;
% \limsup_{n\rightarrow \infty}\frac{\log \log n}{\pen\left(n\right)}=0.
% $$
% Then
% $$
% \left(\widehat{\tau}_{D,0}, \widehat{\tau}_{D,1},\widehat{\tau}_{D,2}\right)=
% \left(\tau_{0}^{\star}, \tau_{1}^{\star}, \tau_{2}^{\star} \right)
% $$
% $\Q^{\star}$ eventually almost surely as $n$ goes to infinity.\\
% Is $\log \log n$ enough ? See van Handel \cite{Han09}.\\

\bibliographystyle{sjs}
\bibliography{references}

\subsection*{Addresses}
\begin{minipage}[t]{6.5cm}
Antonio Galves\\
Universidade de S$\tilde{a}$o Paulo\\
Instituto de Matemática e Estat\'istica - USP\\
Rua do Mat$\tilde{a}$o, 1010\\
CEP 05508-900- S$\tilde{a}$o Paulo\\
Brasil\\
\texttt{galves@usp.br}\\
\end{minipage}\qquad
\begin{minipage}[t]{5cm}
Aur\'elien Garivier\\
CNRS UMR 5141 \\
Telecom-ParisTech\\
LTCI D\'epartement TSI\\ site Dareau\\
46 rue Barrault
75634 Paris cedex 13 \\
France \\
\texttt{garivier@telecom-paristech.fr}\\
\end{minipage}\qquad

\vspace{0.1cm}
\noindent\begin{minipage}[t]{10cm}
Elisabeth Gassiat\\
Laboratoire de {Ma\-th\'e\-ma\-tiques}, CNRS UMR 8628\\
Equipe de Probabilit\'es, Statistique et Mod\'elisation\\
Universit\'e Paris-Sud\\
B\^atiment 425\\
91405 Orsay Cedex\\
France \\
\texttt{elisabeth.gassiat@math.u-psud.fr}\\
\end{minipage}

\vspace{1cm}

\appendix\noindent\textbf{\huge Appendix}

\section{Technical Lemma}\label{sec:tech}
 Let $\P_{U}$ denote the probability distribution of the memoryless source with uniform marginal distribution on $A$. 
 For a context tree $\tau$ and a string $z_{1:k}\in A^k$ denote by $S_{\tau}(\omega,z_{1:k})$ the concatenation of the symbols that are not in context $s$ for any $s\in \tau$, that is
 $S_{\tau}(\omega,z_{1:k})=\bigodot_{i\in I(z_{1:k},\tau)} z_i$.
 Then the  Krichevsky-Trofimov~\cite{k81perf} probability distribution is defined as 
\begin{multline}
\label{modelKT}
\KT_{(\sigma_{0},\sigma_{1},\sigma_{2})}\left(x_{1:n};y_{1:m}\right)=
\P_{U}\left(S_{\sigma_{1}\cup \sigma_{0}}(\omega;x_{1:n})\right)\P_{U}\left(S_{\sigma_{2}\cup \sigma_{0}}(\omega;y_{1:m})\right)
\\
\prod_{s\in\sigma_{0}}
\KT\left(S(s;x_{1:n};y_{1:m})\right)\prod_{s\in\sigma_{1}}
\KT\left(S(s;x_{1:n})\right)\prod_{s\in\sigma_{2}}
\KT\left(S(s;y_{1:m})\right),
\end{multline}
where 
$$
\KT\left(S(s;x_{1:n};y_{1:m})\right)=\frac{\Gamma \left(\frac{|A|}{2}\right)\prod_{a\in A}\Gamma \left(N_{n,x}\left(s,a\right)+N_{n,y}\left(s,a\right)+\frac{1}{2}\right)}{\Gamma \left(\frac{1}{2}\right)^{|A|}\Gamma \left( N_{n,x}\left(s\right)+N_{n,y}\left(s\right)+\frac{|A|}{2}\right)},
$$
$$
\KT\left(S(s;x_{1:n})\right)=\frac{\Gamma \left(\frac{|A|}{2}\right)\prod_{a\in A}\Gamma \left(N_{n,x}\left(s,a\right)+\frac{1}{2}\right)}{\Gamma \left(\frac{1}{2}\right)^{|A|}\Gamma \left( N_{n,x}\left(s\right)+\frac{|A|}{2}\right)},
$$
$$
\KT\left(S(s;y_{1:m})\right)=\frac{\Gamma \left(\frac{|A|}{2}\right)\prod_{a\in A}\Gamma \left(N_{n,y}\left(s,a\right)+\frac{1}{2}\right)}{\Gamma \left(\frac{1}{2}\right)^{|A|}\Gamma \left( N_{n,y}\left(s\right)+\frac{|A|}{2}\right)}.
$$
Recall that for any tree $\sigma$, $D\left(\sigma\right)$ is its depth :
$$
D\left(\sigma\right)=\max \left\{|s|\;:\;s\in \sigma \right\}.
$$
Following Willems \cite{willems1995} (see also \cite{Gassiat2010IT}, and references therein), Jensen's inequality leads to the following result:
\begin{lem}
\label{lemKT}
For any $x_{1:n}$ and any $y_{1:m}$,
\begin{align*}
- \log &\KT_{(\sigma_{0},\sigma_{1},\sigma_{2})}\left(x_{1:n};y_{1:m}\right)  \leq  - \ell_{n,m}\left(\sigma_{0},\sigma_{1},\sigma_{2}\right) \\
&+ \left[D\left(\sigma_{0}\cup \sigma_{1}\right)+D\left(\sigma_{0}\cup \sigma_{2}\right)+ |\sigma_{0}| + |\sigma_{1}| + |\sigma_{2}|\right]\log |A|\\
& + \frac{|A|-1}{2}\left\{|\sigma_{0}| \log \left(\frac{n+m}{|\sigma_{0}|}\right)+|\sigma_{1}| \log \left(\frac{n}{|\sigma_{1}|}\right)+|\sigma_{2}| \log \left(\frac{m}{|\sigma_{2}|}\right)
\right\}
\end{align*}
\end{lem}

\section{Proof of Theorem~\ref{th:consi}}
\label{sec:proof}
The proof is divided into four parts.
\begin{enumerate}
\item
We first prove that eventually almost surely, 
$|\widehat{\sigma}_{0}| \leq k_{n}$ and $|\widehat{\sigma}_{1}| \leq k_{n}$ and $|\widehat{\sigma}_{2}| \leq k_{n}$ with
$$
k_{n}= \frac{\log n}{\log \log \log n}.
$$
For any $(\sigma_{0},\sigma_{1},\sigma_{2})$ satisfying (\ref{T1}), (\ref{T2}) and (\ref{T3}), define
$B_{(\sigma_{0},\sigma_{1},\sigma_{2})}$ as the set of $(x_{1:n},y_{1:m})$ in $A^{n+m}$ such that 
%$$\Q^{\star}\left((\widehat{\sigma}_{0}, \widehat{\sigma}_{1},\widehat{\sigma}_{2})=(\sigma_{0},\sigma_{1},\sigma_{2})\right) =\Q^{\star}\left((X_{1:n},Y_{1:m})\in B_{(\sigma_{0},\sigma_{1},\sigma_{2})}\right).$$
\[(X_{1:n},Y_{1:m})=(x_{1:n},y_{1:m}) \Leftrightarrow (\widehat{\sigma}_{0}, \widehat{\sigma}_{1},\widehat{\sigma}_{2})=(\sigma_{0},\sigma_{1},\sigma_{2}),\]
so that
\begin{multline*}
\Q^{\star}\left((\widehat{\sigma}_{0}, \widehat{\sigma}_{1},\widehat{\sigma}_{2})=(\sigma_{0},\sigma_{1},\sigma_{2})\right)\\
=\sum_{(x_{1:n},y_{1:m})\in B_{(\sigma_{0},\sigma_{1},\sigma_{2})}}\Q^{\star}\left((X_{1:n},Y_{1:m})=(x_{1:n},y_{1:m})\right).
\end{multline*}

If $(X_{1:n},Y_{1:m})\in B_{(\sigma_{0},\sigma_{1},\sigma_{2})}$, then
\begin{multline*}
\ell_{n,m}\left(\sigma_{0},\sigma_{1},\sigma_{2}\right)-\frac{(|A|-1)}{2}(|\sigma_{0}|\pen(n+m)+|\sigma_{1}|\pen(n)+|\sigma_{2}|\pen(m))\\
\geq
\ell_{n,m}\left(\sigma_{0}^{\star},\sigma_{1}^{\star},\sigma_{2}^{\star}\right)-\frac{(|A|-1)}{2}(|\sigma_{0}^{\star}|\pen(n+m)+|\sigma_{1}^{\star}|\pen(n)+|\sigma_{2}^{\star}|\pen(m)),
\end{multline*}
%and since
%$$
%\ell_{n,m}\left(\sigma_{0}^{\star},\sigma_{1}^{\star},\sigma_{2}^{\star}\right) \geq \log \Q^{\star}\left((x_{1:n},y_{1:m})\right)
%$$
and using Lemma \ref{lemKT}, if $(x_{1:n},y_{1:m})\in B_{(\sigma_{0},\sigma_{1},\sigma_{2})}$, then
\begin{align*}
\Q^{\star}&\left((x_{1:n},y_{1:m})\right)  \leq  2^{\ell_{n,m}\left(\sigma_{0}^{\star},\sigma_{1}^{\star},\sigma_{2}^{\star}\right)}
\\ &\leq 2^{\ell_{n,m}\left(\sigma_{0},\sigma_{1},\sigma_{2}\right)+\frac{(|A|-1)}{2}\left(
\left(|\sigma_{0}^*| - t_0\right)\pen(n+m)+\left(|\sigma_{1}^*|-t_1\right)\pen(n)+\left(|\sigma_{2}^*|-t_2\right)\pen(m)
\right)}
\\ &\leq \KT_{(\sigma_{0},\sigma_{1},\sigma_{2})}\left(x_{1:n};y_{1:m}
\right)
2^{H(n,m,t_{0},t_{1},t_{2})}
\end{align*}
with $t_{i}=|\sigma_{i}|$, $i=0,1,2$, and
\begin{align*}
H\big(&n,m,t_{0},t_{1},t_{2}\big)=\\
&\frac{|A|-1}{2}\left\{t_{0} \log \left(\frac{n+m}{t_{0}}\right)+t_{1} \log \left(\frac{n}{t_{1}}\right)+t_{2}\log \left(\frac{m}{t_{2}}\right)
\right\}\\
&+\frac{(|A|-1)}{2}\left((|\sigma_{0}^{\star}|-t_{0})\pen(n+m)+(|\sigma_{1}^{\star}|-t_{1})\pen(n)+(|\sigma_{2}^{\star}|-t_{2})\pen(m)\right)\\
&+\left[3 t_{0}+2 t_{1}+2 t_{2}\right]\log |A|\\
&=\frac{|A|-1}{2}\Big\{-t_{0} \log t_{0}-t_{1} \log t_{1}- t_{2}\log t_{2}+|\sigma_{0}^{\star}|\log \left(n+m\right)+ \\ 
&|\sigma_{1}^{\star}|\log \left(n\right)+|\sigma_{2}^{\star}|\log \left(m\right)\big\}
+\left[3 t_{0}+2 t_{1}+2 t_{2}\right]\log |A|
\end{align*}
using $\pen (\cdot)=\log (\cdot)$ and using that for a complete tree $\sigma$, $D(\sigma ) \leq |\sigma|$. \\
Thus,
$$
\Q^{\star}\left((\widehat{\sigma}_{0}, \widehat{\sigma}_{1},\widehat{\sigma}_{2})=(\sigma_{0},\sigma_{1},\sigma_{2})\right)
\leq
2^{H(n,m,t_{0},t_{1},t_{2})},
$$
and
\begin{multline*}
\Q^{\star}\left(|\widehat{\sigma}_{0}|\geq k_{n} \text{ or } |\widehat{\sigma}_{1}|\geq k_{n} \text{ or } |\widehat{\sigma}_{2}|\geq k_{n} \right)\\
\leq \sum_{t_{0}=k_{n}+1}^{n\vee m} \sum_{t_{1},t_{2}=0}^{n\vee m} F\left(t_{0},t_{1},t_{2}\right)2^{H(n,m,t_{0},t_{1},t_{2})}\\
+ \sum_{t_{1}=k_{n}+1}^{n\vee m} \sum_{t_{0},t_{2}=0}^{n\vee m} F\left(t_{0},t_{1},t_{2}\right)2^{H(n,m,t_{0},t_{1},t_{2})}\\
+ \sum_{t_{2}=k_{n}+1}^{n\vee m} \sum_{t_{0},t_{1}=0}^{n\vee m} F\left(t_{0},t_{1},t_{2}\right)2^{H(n,m,t_{0},t_{1},t_{2})}
\end{multline*}
where $F\left(t_{0},t_{1},t_{2}\right)$ is the number of $(\sigma_{0},\sigma_{1},\sigma_{2})$ satisfying (\ref{T1}), (\ref{T2}) and (\ref{T3}) and such that
$|\sigma_{0}|=t_{0}$, $|\sigma_{1}|=t_{1}$, and $|\sigma_{2}|=t_{2}$.\\
But the number of complete trees with $t$ elements is upper bounded by $16^{t}$, see \cite{garivier2006}, so that, denoting by $\binom{b}{a}\leq 2^b$ the binomial coefficient, one has
\begin{eqnarray*}
F\left(t_{0},t_{1},t_{2}\right) &\leq & \binom{t_{0}+t_{1}}{t_{0}} 16^{t_{0}+t_{1}}\binom{t_{0}+t_{2}}{t_{0}} 16^{t_{0}+t_{2}}\\
&\leq & 16^{4t_{0}+2t_{1}+2t_{2}}.
\end{eqnarray*}
Using the fact that for any constant $a$, $-t\log t + a t$ is bounded on $\R^{+}$, and using (\ref{equilibre}) one gets that for some constants
$C_1$, $C_2$ and $C_3$,
$$
\Q^{\star}\left(|\widehat{\sigma}_{0}|\geq k_{n} \text{ or } |\widehat{\sigma}_{1}|\geq k_{n} \text{ or } |\widehat{\sigma}_{2}|\geq k_{n} \right)
\leq C_1 2^{-C_2 k_{n}\log k_{n} + C_3 \log n}.
$$
But 
$$
\lim_{n\rightarrow +\infty}\frac{k_{n}\log k_{n}}{\log n}=+\infty
$$
so that one gets that for another constant $C$,
$$
\Q^{\star}\left(|\widehat{\sigma}_{0}|\geq k_{n} \text{ or } |\widehat{\sigma}_{1}|\geq k_{n} \text{ or } |\widehat{\sigma}_{2}|\geq k_{n} \right)
\leq \frac{C}{n^{2}}
$$
and using Borel-Cantelli's Lemma, we obtain that
$\Q^{\star}$-eventually almost surely, 
$|\widehat{\sigma}_{0}| \leq k_{n}$ and $|\widehat{\sigma}_{1}| \leq k_{n}$ and $|\widehat{\sigma}_{2}| \leq k_{n}$.

\item
We prove that $\Q^{\star}$-eventually almost surely, no context is overestimated.\\
It is sufficient to prove that, $\Q^{\star}$-almost surely, if $(\sigma_{0},\sigma_{1},\sigma_{2})$ satisfy (\ref{T1}), (\ref{T2}) and (\ref{T3}) and are such that for some $i$, $\sigma_{i}$ contains some string that has a proper suffix in $\sigma_{i}^{\star}$, there exists $(\bar{\sigma}_{0},\bar{\sigma}_{1},\bar{\sigma}_{2})$ satisfying (\ref{T1}), (\ref{T2}) and (\ref{T3})
and such that, eventually,  $C_{n,m}(\bar{\sigma}_{0},\bar{\sigma}_{1},\bar{\sigma}_{2})>C_{n,m}(\sigma_{0},\sigma_{1},\sigma_{2})$, so that
$\left(\widehat{\sigma}_{0}, \widehat{\sigma}_{1},\widehat{\sigma}_{2}\right)\neq (\sigma_{0},\sigma_{1},\sigma_{2})$.\\
Consider first the case where $\sigma_{0}^{\star}$ is overestimated. Let  $(\sigma_{0},\sigma_{1},\sigma_{2})$ satisfy (\ref{T1}), (\ref{T2}) and (\ref{T3}) and be such that $\sigma_{0}$ contains some string that has a proper suffix in $\sigma_{0}^{\star}$. Let $s=av$, $a\in A$, be the longest such string, and let $u\in\sigma_{0}^{\star}$ be the corresponding suffix of $v$. For $i\in\{0,1,2\}$, let
$S_{i}=A^{+}v\cap \sigma_{i}$ and %,\;L_{i}=\left\vert S_{i} \right\vert.$
define 
\[
\bar{\sigma}_{0}=\left(\sigma_{0}\backslash S_{0}\right)\cup \{v\}\;,\qquad
\bar{\sigma}_{1}=\left(\sigma_{1}\backslash S_{1}\right)\;,\qquad
\bar{\sigma}_{2}=\left(\sigma_{2}\backslash S_{2}\right)\;.
\]
Then
\begin{align*}
&C_{n,m}(\bar{\sigma}_{0},\bar{\sigma}_{1},\bar{\sigma}_{2})-C_{n,m}(\sigma_{0},\sigma_{1},\sigma_{2})\\
&\hspace{-0.3cm}=\quad
\sum_{b\in A} \left[N_{n,X}\left(v,b\right)+N_{m,Y}\left(v,b\right)\right]\log\left(\frac{N_{n,X}\left(v,b\right)+N_{m,Y}\left(v,b\right)}{N_{n,X}\left(v\right)+N_{m,Y}\left(v\right)}\right)\\
& \qquad - \frac{|A|-1}{2}\log \left(n+m\right)\\
&
-\sum_{w\in S_{0}}\Bigg\{\sum_{b\in A} \left[N_{n,X}\left(w,b\right)+N_{m,Y}\left(w,b\right)\right]\log\left(\frac{N_{n,X}\left(w,b\right)+N_{m,Y}\left(w,b\right)}{N_{n,X}\left(w\right)+N_{m,Y}\left(w\right)}\right)\\
&\qquad -\frac{|A|-1}{2}\log \left(n+m\right)\Bigg\}\\
&
-\sum_{w\in S_{1}}\left\{\sum_{b\in A} N_{n,X}\left(w,b\right)\log\left(\frac{N_{n,X}\left(w,b\right)}{N_{n,X}\left(w\right)}\right)- \frac{|A|-1}{2}\log \left(n\right)\right\}\\
&
-\sum_{w\in S_{2}}\left\{\sum_{b\in A} N_{m,Y}\left(w,b\right)\log\left(\frac{N_{m,Y}\left(w,b\right)}{N_{m,Y}\left(w\right)}\right)- \frac{|A|-1}{2}\log \left(m\right)\right\}
\end{align*}
By definition of the maximum likelihood, the above expression is lower-bounded by:
\begin{align*}
&C_{n,m}(\bar{\sigma}_{0},\bar{\sigma}_{1},\bar{\sigma}_{2})-C_{n,m}(\sigma_{0},\sigma_{1},\sigma_{2})\\
&\hspace{-0.3cm}\geq \quad
\sum_{b\in A} \left[N_{n,X}\left(v,b\right)+N_{m,Y}\left(v,b\right)\right]\log\left(Q^{\star}_{X}\left(b \vert v \right)\right) - \frac{|A|-1}{2}\log \left(n+m\right)\\
& 
-\sum_{w\in S_{0}}\Bigg\{\sum_{b\in A} \left[N_{n,X}\left(w,b\right)+N_{m,Y}\left(w,b\right)\right]\log\left(\widehat{Q}_{XY}\left(b \vert w \right)\right)\\
&\qquad- \frac{|A|-1}{2}\log \left(n+m\right)\Bigg\}\\
&
-\sum_{w\in S_{1}}\left\{\sum_{b\in A} N_{n,X}\left(w,b\right)\log\left(\widehat{Q}_{X}\left(b \vert w \right)\right)- \frac{|A|-1}{2}\log \left(n\right)\right\}\\
&
-\sum_{w\in S_{2}}\left\{\sum_{b\in A} N_{m,Y}\left(w,b\right)\log\left(\widehat{Q}_{Y}\left(b \vert w \right)\right)- \frac{|A|-1}{2}\log \left(m\right)\right\}
\end{align*}
Notice that  \[Q^{\star}_{X}\left(\cdot \vert v \right)=Q^{\star}_{Y}\left(\cdot \vert v \right)=Q^{\star}_{X}\left(\cdot \vert w \right)\] for any $w\in S_{0}\cup S_{1}\cup S_{2}$.\\
It follows from part 1 of the proof that we only need to consider  trees $\sigma_{i}$ such that $|\sigma_{i}|=o(\log n)$. Notice also that since $D(\sigma_{i})=o(\log n)$,
for any $b\in A$,
$$
N_{n,X}\left(v,b\right)=\sum_{w\in S_{0}\cup S_{1}}N_{n,X}\left(w,b\right)+o(\log n),
$$
$$
N_{m,Y}\left(v,b\right)=\sum_{w\in S_{0}\cup S_{2}}N_{m,Y}\left(w,b\right)+o(\log n).
$$
Let $\kl\left(q_1 \vert q_2\right) = \sum_{a\in A}q_1(a)\log\frac{q_1(a)}{q_2(a)}$ denotes the Kullback-Leibler divergence between two probability measures $q_1$ and $q_2$ on $A$, with the convention that $0\log(0/x)=0$ for $x\geq 0$ and $x\log(x/0)=+\infty$ for $x>0$.
Since the minimum of all positive transition probabilities in $\Q^{\star}$ is positive, one gets
\begin{align*}
C_{n,m}&(\bar{\sigma}_{0},\bar{\sigma}_{1},\bar{\sigma}_{2})-C_{n,m}(\sigma_{0},\sigma_{1},\sigma_{2}) \\
&\hspace{-0.3cm}\geq \hspace{0.3cm}
\sum_{w\in S_{0}}\sum_{b\in A} \left[N_{n,X}\left(w,b\right)+N_{m,Y}\left(w,b\right)\right]\log\left(\frac{Q^{\star}_{X}\left(b \vert w \right)}{\widehat{Q}_{XY}\left(b \vert w \right)}\right) \\
&\qquad+ \left(|S_{0}|-1\right) \frac{|A|-1}{2}\log \left(n+m\right)
\\
&
+\sum_{w\in S_{1}}\sum_{b\in A} N_{n,X}\left(w,b\right)\log\left(\frac{Q^{\star}_{X}\left(b \vert w \right)}{\widehat{Q}_{X}\left(b \vert w \right)}\right) +|S_{1}|\frac{|A|-1}{2}\log \left(n\right)\\
&
+\sum_{w\in S_{2}}\sum_{b\in A} N_{m,Y}\left(w,b\right)\log\left(\frac{Q^{\star}_{Y}\left(b \vert w \right)}{\widehat{Q}_{Y}\left(b \vert w \right)}\right)+ |S_{2}|\frac{|A|-1}{2}\log \left(m\right)
\\
&
+ o ( \log n   )\\
& \hspace{-0.5cm}= \hspace{0.5cm}
-\sum_{w\in S_{0}}\left[N_{n,X}\left(w\right)+N_{m,Y}\left(w\right)\right]\kl\left( \widehat{Q}_{XY}\left(\cdot \vert w \right)\vert Q^{\star}_{X}\left(\cdot \vert w \right) \right)\\
&\qquad+\left(|S_{0}|-1\right) \frac{|A|-1}{2}\log \left(n+m\right)
\\
&
-\sum_{w\in S_{1}} N_{n,X}\left(w\right)\kl\left(\widehat{Q}_{X}\left(\cdot \vert w \right) \vert Q^{\star}_{X}\left(\cdot \vert w \right)\right) +|S_{1}|\frac{|A|-1}{2}\log \left(n\right)
\\
&
-\sum_{w\in S_{2}} N_{m,Y}\left(w\right)\kl\left(  \widehat{Q}_{Y}\left(\cdot \vert w \right) \vert Q^{\star}_{Y}\left(\cdot \vert w \right)\right)+ |S_{2}|\frac{|A|-1}{2}\log \left(m\right)\\
&+ o \left( \log n  \right).
\end{align*}
According to typicality Lemma 6.2 of \cite{csiszar2006}, for all $\delta>0$, for all $w$ such that $N_{n,X}(w)\geq 1$ and for all $b\in A$ it holds that, $\Q^{\star}$-eventually almost surely,
\[\left|\widehat{Q}_{X}\left(b \vert w \right) - Q^{\star}_{X}\left(b \vert w \right)\right| \leq \sqrt{\frac{\delta\log(n)}{N_{n,X}\left(w\right)}}\;.\]
Besides, Lemma 6.3 of \cite{csiszar2006} states that 
\[\kl\left( \widehat{Q}_{X}\left(\cdot \vert w \right) \vert Q^{\star }_{X}\left(\cdot \vert w \right)\right) \leq \sum_{b\in A}\frac{\left(\widehat{Q}_{X}\left(b \vert w \right)-Q^{\star }_{X}\left(b \vert w \right)\right)^2}{Q^{\star }_{X}\left(b \vert w \right)}\;. \]
Handling similarly the terms involving $Q^{\star }_{Y}$ and $Q^{\star }_{XY}$, and denoting $q^{\star}_{min}>0$ the minimum of all positive transition probabilities in $\Q^{\star}$, we obtain that for any $\delta >0$, $\Q^{\star}$-eventually almost surely for all possible 
$(\sigma_{0},\sigma_{1},\sigma_{2})$ :
\begin{multline*}
C_{n,m}(\bar{\sigma}_{0},\bar{\sigma}_{1},\bar{\sigma}_{2})-C_{n,m}(\sigma_{0},\sigma_{1},\sigma_{2})\geq \\
-\frac{\delta |A|}{q^{\star}_{min}}|S_{0}| \log \left(n+m\right)+\left(|S_{0}|-1\right) \frac{|A|-1}{2}\log \left(n+m\right)\\
-\frac{\delta |A|}{q^{\star}_{min}}|S_{1}| \log \left(n\right)
+|S_{1}|\frac{|A|-1}{2}\log \left(n\right)\\
-\frac{\delta |A|}{q^{\star}_{min}}|S_{2}| \log \left(m\right)
+ |S_{2}|\frac{|A|-1}{2}\log \left(m\right)
\end{multline*}
which is positive, for all possible 
$(\sigma_{0},\sigma_{1},\sigma_{2})$, $\Q^{\star}$-eventually almost surely. This follows from the fact that we consider complete context trees, and therefore $|S_{0}|\geq 1$, $|S_{0}|+|S_{1}| \geq |A|$ and $|S_{0}|+|S_{2}| \geq |A|$.\\
Consider now the case where $\sigma_{i}^{\star}$, $i=1$ or $i=2$ is overestimated.
Let  $(\sigma_{0},\sigma_{1},\sigma_{2})$ satisfy (\ref{T1}), (\ref{T2}) and (\ref{T3}) and be such that $\sigma_{i}$ contains some string that has a proper suffix in $\sigma_{i}^{\star}$. Let $s=av$, $a\in A$, be the longest such string, and let $u\in\sigma_{i}^{\star}$ be the corresponding suffix of $v$.
For $i=0,1,2$, let again, $S_{i}=A^{+}v\cap \sigma_{i}$. %,\;|S_{i}|=\left\vert S_{i} \right\vert.$$
Then, either $S_{0}=\emptyset$, and the problem boils down the the overestimation of a single tree: the consistency result of \cite{csiszar2006} applies and shows that denoting
\[
\bar{\sigma}_{i}=\left(\sigma_{1}\backslash S_{i}\right)\cup \left\{v\right\} \;,\quad
\bar{\sigma}_{j}=\sigma_{j}, j\neq i\;, 
\]
we have $C_{n,m}(\bar{\sigma}_{0},\bar{\sigma}_{1},\bar{\sigma}_{2})>C_{n,m}(\sigma_{0},\sigma_{1},\sigma_{2})$  $\Q^{\star}$-eventually almost surely.
Or $\sigma_{0}^{\star}$ has also been overestimated, so that one may apply the previous proof.
\item
Consider now the underestimation case. 
If $\sigma_{0}$ has been underestimated, there exists $s\in \sigma_{0}$ which is a proper suffix of $s_{0}\in\sigma_{0}^{\star}$. For $i=0,1,2$, let $
S_{i}=A^{+}s\cap \sigma_{i}^{\star}$, and define%,\;|S_{i}|=\left\vert S_{i} \right\vert.$$
\[
\bar{\sigma}_{0}=\left(\sigma_{0}\backslash \{s\}\right)\cup S_{0}\;,\quad
\bar{\sigma}_{1}=\sigma_{1}\cup S_{1}\;,\quad
\bar{\sigma}_{2}=\sigma_{2}\cup S_{2}.
\]
Then
\begin{align*}
C_{n,m}(&\bar{\sigma}_{0},\bar{\sigma}_{1},\bar{\sigma}_{2})-C_{n,m}(\sigma_{0},\sigma_{1},\sigma_{2})\\
=&
\sum_{w\in S_{0}}\Bigg\{\sum_{b\in A} \left[N_{n,X}\left(w,b\right)+N_{m,Y}\left(w,b\right)\right]\log\left(\frac{N_{n,X}\left(w,b\right)+N_{m,Y}\left(w,b\right)}{N_{n,X}\left(w\right)+N_{m,Y}\left(w\right)}\right)\\
&\qquad- \frac{|A|-1}{2}\log \left(n+m\right)\Bigg\}\\
&
+\sum_{w\in S_{1}}\left\{\sum_{b\in A} N_{n,X}\left(w,b\right)\log\left(\frac{N_{n,X}\left(w,b\right)}{N_{n,X}\left(w\right)}\right)- \frac{|A|-1}{2}\log \left(n\right)\right\}\\
&
+\sum_{w\in S_{2}}\left\{\sum_{b\in A} N_{m,Y}\left(w,b\right)\log\left(\frac{N_{m,Y}\left(w,b\right)}{N_{m,Y}\left(w\right)}\right)- \frac{|A|-1}{2}\log \left(m\right)\right\}\\
&
-\sum_{b\in A} \left[N_{n,X}\left(s,b\right)+N_{m,Y}\left(s,b\right)\right]\log\left(\frac{N_{n,X}\left(s,b\right)+N_{m,Y}\left(s,b\right)}{N_{n,X}\left(s\right)+N_{m,Y}\left(v\right)}\right) \\
&\qquad+ \frac{|A|-1}{2}\log \left(n+m\right)
\end{align*}
Notice that for any string $u$, for any $b\in A$, 
$
\frac{1}{n} N_{n,X}\left(u,b\right)
$ and $
\frac{1}{n} N_{n,X}\left(u\right)
$
converge $\Q^{\star}$ almost surely to $Q^{\star}_{X}\left(ub\right)$ and $Q^{\star}_{X}\left(u\right)$ respectively, and
$
\frac{1}{n} N_{m,Y}\left(u,b\right)
$ and $
\frac{1}{n} N_{m,Y}\left(u \right)
$
converge $\Q^{\star}$ almost surely to $\frac{1}{c}Q^{\star}_{Y}\left(ub\right)$ and $\frac{1}{c}Q^{\star}_{Y}\left(u\right)$, respectively.

Thus, $\Q^{\star}$ almost surely,
\begin{align*}
&C_{n,m}(\bar{\sigma}_{0},\bar{\sigma}_{1},\bar{\sigma}_{2})-C_{n,m}(\sigma_{0},\sigma_{1},\sigma_{2})=- O\left(\log n \right)\\
&
+ n 
\sum_{w\in S_{0}}\sum_{b\in A} \left[Q^{\star}_{X}\left(wb\right)+\frac{1}{c}Q^{\star}_{Y}\left(wb\right)\right]\log\left(\frac{Q^{\star}_{X}\left(wb\right)+\frac{1}{c}Q^{\star}_{Y}\left(wb\right)}{Q^{\star}_{X}\left(w\right)+\frac{1}{c}Q^{\star}_{Y}\left(w\right)}\right)\\
&
+n \sum_{w\in S_{1}}\sum_{b\in A} Q^{\star}_{X}\left(wb\right)\log\left(\frac{Q^{\star}_{X}\left(wb\right)}{Q^{\star}_{X}\left(w\right)}\right)\\
&
+n \sum_{w\in S_{2}}\sum_{b\in A} \frac{1}{c}Q^{\star}_{Y}\left(wb\right)\log\left(\frac{Q^{\star}_{Y}\left(wb\right)}{Q^{\star}_{Y}\left(w\right)}\right)\\
&
-n\sum_{b\in A} \left[Q^{\star}_{X}\left(sb\right)+\frac{1}{c}Q^{\star}_{Y}\left(sb\right)\right]\log\left(\frac{Q^{\star}_{X}\left(sb\right)+\frac{1}{c}Q^{\star}_{Y}\left(sb\right)}{Q^{\star}_{X}\left(s\right)+\frac{1}{c}Q^{\star}_{Y}\left(s\right)}\right) + o\left(n\right)
\\
& =- O\left(\log n \right)+ o\left(n\right)
+n \sum_{w\in S_{0}\cup S_{1}}\sum_{b\in A} Q^{\star}_{X}\left(wb\right)\log\left(\frac{Q^{\star}_{X}\left(wb\right)}{Q^{\star}_{X}\left(w\right)}\right)\\
&
+n \sum_{w\in S_{0}\cup S_{2}}\sum_{b\in A} \frac{1}{c}Q^{\star}_{Y}\left(wb\right)\log\left(\frac{Q^{\star}_{Y}\left(wb\right)}{Q^{\star}_{Y}\left(w\right)}\right)\\
&
-n\sum_{b\in A} \left[Q^{\star}_{X}\left(sb\right)+\frac{1}{c}Q^{\star}_{Y}\left(sb\right)\right]\log\left(\frac{Q^{\star}_{X}\left(sb\right)+\frac{1}{c}Q^{\star}_{Y}\left(sb\right)}{Q^{\star}_{X}\left(s\right)+\frac{1}{c}Q^{\star}_{Y}\left(s\right)}\right) 
\end{align*}
because  for $w\in S_{0}$, $Q^{\star}_{X}\left(wb\right)=Q^{\star}_{Y}\left(wb\right)$.
Since
$$
\sum_{w\in S_{0}\cup S_{1}}Q^{\star}_{X}\left(w\right)=Q^{\star}_{X}\left(s\right),
$$
for any $b\in A$, Jensen's inequality  implies that %on the convex mapping $x\to \log(1/x)$ with weights $Q_X^{\star}(wb)/X_X^{\star}(sb)$ leads to
$$
\sum_{w\in S_{0}\cup S_{1}} Q^{\star}_{X}\left(wb\right)\log\left(\frac{Q^{\star}_{X}\left(wb\right)}{Q^{\star}_{X}\left(w\right)}\right)
\geq Q^{\star}_{X}\left(sb\right)\log\left(\frac{Q^{\star}_{X}\left(sb\right)}{Q^{\star}_{X}\left(s\right)}\right),
$$
and the inequality is strict for at least one $b\in A$, for otherwise, $s$ would be a context for $\Q^{\star}_{X}$.  
Similarly for any $b\in A$,
$$
\sum_{w\in S_{0}\cup S_{2}} Q^{\star}_{Y}\left(wb\right)\log\left(\frac{Q^{\star}_{Y}\left(wb\right)}{Q^{\star}_{Y}\left(w\right)}\right)
\geq Q^{\star}_{Y}\left(sb\right)\log\left(\frac{Q^{\star}_{Y}\left(sb\right)}{Q^{\star}_{Y}\left(s\right)}\right).
$$
Using the concavity of the entropy function
\begin{align*}
&\sum_{b\in A}Q^{\star}_{X}\left(sb\right)\log\left(\frac{Q^{\star}_{X}\left(sb\right)}{Q^{\star}_{X}\left(s\right)}\right) 
+ \frac{1}{c} \sum_{b\in A}Q^{\star}_{Y}\left(sb\right)\log\left(\frac{Q^{\star}_{Y}\left(sb\right)}{Q^{\star}_{Y}\left(s\right)}\right)\\
&\geq
\sum_{b\in A}\left(Q^{\star}_{X}\left(sb\right)+ \frac{1}{c} Q^{\star}_{Y}\left(sb\right)\right)\log\left(\frac{Q^{\star}_{X}\left(sb\right)+ \frac{1}{c} Q^{\star}_{Y}\left(sb\right)}{Q^{\star}_{X}\left(s\right)+ \frac{1}{c} Q^{\star}_{Y}\left(s\right)}\right),
\end{align*}
so that there exists $\delta >0$ such that 
$$
C_{n,m}(\bar{\sigma}_{0},\bar{\sigma}_{1},\bar{\sigma}_{2})-C_{n,m}(\sigma_{0},\sigma_{1},\sigma_{2}) \geq n\delta
$$
$\Q^{\star}$-eventually almost surely.
\\
If $\sigma_{i}$, $i=1$ or $i=2$  has been underestimated, then the problem boils down to the standard underestimation of a single context tree.
Defining (with obvious notation) 
 \[
\bar{\sigma}_{i}=\left(\sigma_{1}\backslash \left\{s\right\}\right)\cup S_{i}\cup S_{0}\;,\quad
\bar{\sigma}_{j}=\sigma_{j}, j\neq i\;,
\]
it is proved in~\cite{csiszar2006}, Section III, that  $\Q^{\star}$-eventually almost surely,
 $C_{n,m}(\bar{\sigma}_{0},\bar{\sigma}_{1},\bar{\sigma}_{2})>C_{n,m}(\sigma_{0},\sigma_{1},\sigma_{2})$.

\item
We have thus proved that, for $i=1$ and $i=2$, $\widehat{\sigma}_{0}\cup \widehat{\sigma}_{i}={\sigma}^{\star}_{0}\cup {\sigma}^{\star}_{i}$, $\Q^{\star}$-eventually almost surely. Let $(\sigma_{0},\sigma_{1},\sigma_{2})$  satisfy (\ref{T1}), (\ref{T2}) and (\ref{T3}) and be such that, for $i=1$ and $i=2$, $\sigma_{0}\cup \sigma_{i}={\sigma}^{\star}_{0}\cup {\sigma}^{\star}_{i}$. There remains to check that $\Q^{\star}$ almost surely, if there exists a string $s$
such that
\begin{itemize}
\item
$s\in \sigma_{0}$, but $s\in \sigma_{1}^{\star}$ and $s\in \sigma_{2}^{\star}$,
\item
or  $s\in \sigma_{1}$ and $s\in \sigma_{2}$, but  $s\in \sigma_{0}^{\star}$,
\end{itemize}
then $(\widehat{\sigma_{0}},\widehat{\sigma_{1}},\widehat{\sigma_{2}})\neq (\sigma_{0},\sigma_{1},\sigma_{2})$ eventually.\\ 
Consider first the case where $s\in \sigma_{0}$, but $s\in \sigma_{1}^{\star}$ and $s\in \sigma_{2}^{\star}$. 
Define 
\[\bar{\sigma}_{0}=\left(\sigma_{0}\backslash \{s\}\right)\;,\quad
\bar{\sigma}_{1}=\sigma_{1}\cup \{s\}\;,\quad
\bar{\sigma}_{2}=\sigma_{2}\cup \{s\}\;.
\]
Then
\begin{align*}
&C_{n,m}(\bar{\sigma}_{0},\bar{\sigma}_{1},\bar{\sigma}_{2})-C_{n,m}(\sigma_{0},\sigma_{1},\sigma_{2})=\\
&
+\sum_{b\in A} N_{n,X}\left(s,b\right)\log\left(\frac{N_{n,X}\left(s,b\right)}{N_{n,X}\left(s\right)}\right)\\
&
\sum_{b\in A} N_{m,Y}\left(s,b\right)\log\left(\frac{N_{m,Y}\left(s,b\right)}{N_{m,Y}\left(s\right)}\right)\\
&
-\sum_{b\in A} \left[N_{n,X}\left(s,b\right)+N_{m,Y}\left(s,b\right)\right]\log\left(\frac{N_{n,X}\left(s,b\right)+N_{m,Y}\left(s,b\right)}{N_{n,X}\left(s\right)+N_{m,Y}\left(s\right)}\right) \\
&+ \frac{|A|-1}{2}\left\{\log \left(n+m\right)-\log n -\log m\right\}\\
&= n\left\{\sum_{b\in A}Q^{\star}_{X}\left(sb\right)\log\left(\frac{Q^{\star}_{X}\left(sb\right)}{Q^{\star}_{X}\left(s\right)}\right) 
+ \frac{1}{c} \sum_{b\in A}Q^{\star}_{Y}\left(sb\right)\log\left(\frac{Q^{\star}_{Y}\left(sb\right)}{Q^{\star}_{Y}\left(s\right)}\right)\right.\\
&-\left.
\sum_{b\in A}\left(Q^{\star}_{X}\left(sb\right)+ \frac{1}{c} Q^{\star}_{Y}\left(sb\right)\right)\log\left(\frac{Q^{\star}_{X}\left(sb\right)+ \frac{1}{c} Q^{\star}_{Y}\left(sb\right)}{Q^{\star}_{X}\left(s\right)+ \frac{1}{c} Q^{\star}_{Y}\left(s\right)}\right)+o(1)\right\}\\
&-O\left(\log n\right)
\end{align*}
$\Q^{\star}$ almost surely. But the quantity into brackets is positive by the strict concavity of the entropy function, unless for any $b\in A$,
$Q^{\star}_{X}(b\vert s )=Q^{\star}_{Y}(b\vert s )$ which would mean that $s\in \sigma_{0}^{\star}$.\\
Consider now the case where  $s\in \sigma_{1}$ and $s\in \sigma_{2}$, but  $s\in \sigma_{0}^{\star}$.
Define 
\begin{align*}
\bar{\sigma}_{0}&=\sigma_{0}\cup \{s\},\\
\bar{\sigma}_{1}&=\left(\sigma_{1}\backslash \{s\}\right),\\
\bar{\sigma}_{2}&=\left(\sigma_{2}\backslash \{s\}\right).
\end{align*}
\begin{align*}
&C_{n,m}(\bar{\sigma}_{0},\bar{\sigma}_{1},\bar{\sigma}_{2})-C_{n,m}(\sigma_{0},\sigma_{1},\sigma_{2})=\\
&
\sum_{b\in A} \left[N_{n,X}\left(s,b\right)+N_{m,Y}\left(s,b\right)\right]\log\left(\frac{N_{n,X}\left(s,b\right)+N_{m,Y}\left(s,b\right)}{N_{n,X}\left(s\right)+N_{m,Y}\left(s\right)}\right) \\
&-\sum_{b\in A} N_{n,X}\left(s,b\right)\log\left(\frac{N_{n,X}\left(s,b\right)}{N_{n,X}\left(s\right)}\right)\\
&
-\sum_{b\in A} N_{m,Y}\left(s,b\right)\log\left(\frac{N_{m,Y}\left(s,b\right)}{N_{m,Y}\left(s\right)}\right)\\
&+ \frac{|A|-1}{2}\left\{\log n +\log m - \log \left(n+m\right)\right\}.
\end{align*}
Using Taylor expansion until second order of $u \log u$, one gets
\begin{align*}
C_{n,m}&(\bar{\sigma}_{0},\bar{\sigma}_{1},\bar{\sigma}_{2})-C_{n,m}(\sigma_{0},\sigma_{1},\sigma_{2})\\
&\hspace{-0.5cm}= \hspace{0.1cm}\Bigg\{
\frac{1}{2}\sum_{b\in A} \frac{(\left[N_{n,X}\left(s,b\right)+N_{m,Y}\left(s,b\right)\right]-\left[N_{n,X}\left(s\right)+N_{m,Y}\left(s\right)\right] Q^{\star}_{X}(b \vert s ))^{2}}{\left[N_{n,X}\left(s\right)+N_{m,Y}\left(s\right)\right] Q^{\star}_{X}(b \vert s )}\\
& -\frac{1}{2}\sum_{b\in A} \frac{(N_{n,X}\left(s,b\right)-N_{n,X}\left(s\right) Q^{\star}_{X}(b \vert s ))^{2}}{N_{n,X}\left(s\right) Q^{\star}_{X}(b \vert s )}\\
&-\frac{1}{2}\sum_{b\in A} \frac{(N_{m,Y}\left(s,b\right)-N_{m,Y}\left(s\right) Q^{\star}_{Y}(b \vert s ))^{2}}{N_{m,Y}\left(s\right) Q^{\star}_{Y}(b \vert s )}
\Bigg\}\left(1+o(1)\right)\\
&+ \frac{|A|-1}{2}\left\{\log n +\log m - \log \left(n+m\right)\right\}.
\end{align*}
The sequences 
\begin{align*}
 &\left(N_{n,X}\left(s,b\right)-N_{n,X}\left(s\right) Q^{\star}_{X}(b \vert s )\right)_{n\geq 0}\;,\\ 
 &\left(N_{m,Y}\left(s,b\right)-N_{m,Y}\left(s\right) Q^{\star}_{Y}(b \vert s )\right)_{m\geq 0} \;,
 \end{align*}
 are martingales with respect to the natural filtration.
Thus, it follows from the law of iterated logarithm for martingales~\cite{neveu72} that, $\Q^{\star}$ almost surely,
\begin{multline*}
C_{n,m}(\bar{\sigma}_{0},\bar{\sigma}_{1},\bar{\sigma}_{2})-C_{n,m}(\sigma_{0},\sigma_{1},\sigma_{2})=O\left(\log \log n\right)\\+
 \frac{|A|-1}{2}\left\{\log n +\log m - \log \left(n+m\right)\right\},
\end{multline*}
so that $\Q^{\star}$ almost surely,
$$
C_{n,m}(\bar{\sigma}_{0},\bar{\sigma}_{1},\bar{\sigma}_{2})-C_{n,m}(\sigma_{0},\sigma_{1},\sigma_{2}) >0
$$
eventually.
This ends the proof of Theorem \ref{th:consi}.
\end{enumerate}

\section*{Acknowledgments}
This work is part of USP project MaCLinC, ``Mathematics, computation,
language and the brain", USP/COFECUB project
``Stochastic systems with interactions of variable range'' and CNPq project
\emph{Rhythmic patterns, prosodic domains and
  probabilistic modeling in Portuguese Corpora} (grant number
485999/2007-2). This paper was partially supported by CAPES grant AUXPE-PAE-598/2011. A.~Galves is partially supported by a CNPq fellowship (grant
305447/2008-4). The authors thank P. Weyer-Brown for his help in English.

\end{document}